\documentclass[a4paper,10pt]{amsart}

\usepackage{amstext}
\usepackage{amsmath}
\usepackage{ifthen}
\usepackage{amscd}
\usepackage{amssymb}
\usepackage[english]{babel}
\usepackage[T1]{fontenc}
\usepackage[utf8]{inputenc}
\usepackage[all]{xy}
\usepackage{graphicx}
\usepackage{enumerate}
\usepackage{xspace}
\usepackage{pdfsync}
\usepackage{epic}

\newenvironment{demo}[1][]{\ifthenelse{\equal{#1}{}}{\noindent\textbf{Proof:}\xspace}{\noindent\textbf{Proof #1:}\xspace}}{$\square$}

\newtheoremstyle{break}
  {}
  {}
  {\itshape}
  {}
  {\bfseries}
  {.}
  {3pt}
  {}
  
\newtheoremstyle{rq}
  {}
  {}
  {\slshape}
  {}
  {\bfseries}
  {.}
  {3pt}
  {}

\newtheoremstyle{exemple}
  {}
  {}
  {\upshape}
  {}
  {\bfseries}
  {.}
  {3pt}
  {}

\newtheoremstyle{fact}
  {}
  {}
  {\slshape}
  {}
  {\bfseries}
  {.}
  {2pt}
  {}

\theoremstyle{fact}
\newtheorem*{fact*}{Fact}
\theoremstyle{break}
\newtheorem{thm}{Theorem}[section]
\newtheorem{conj}{Conjecture}[section]
\newtheorem{cor}{Corollary}[section]
\newtheorem{lem}{Lemma}[section]
\newtheorem{prop}{Proposition}[section]
\newtheorem{defi}{Definition}[section]
\theoremstyle{rq}
\newtheorem{rem}{Remark}[section]

\theoremstyle{exemple}
\newtheorem{ex}{Example}[section]

\newcommand{\RR}{\mathbb R}

\newcommand{\CC}{\mathbb C}
\newcommand{\QQ}{\mathbb Q}
\newcommand{\ZZ}{\mathbb Z}

\newcommand{\PP}{\mathbb P}

\newcommand{\XX}{\ensuremath{\mathcal{X}}\xspace}

\newcommand{\lie}{\mathcal{L}}
\newcommand{\kod}[1]{\kappa(#1)}
\newcommand{\gf}[1]{\pi_1(#1)}

\newcommand{\To}{\longrightarrow}
\newcommand{\abs}[1]{\left\vert#1\right\vert}

\newcommand{\set}[1]{\left\{#1\right\}}

\newcommand{\dimm}[1]{\mathrm{dim}(#1)}

\newcommand{\wtx}{\tilde{X}}

\newcommand{\merom}[3]{\ensuremath{#1:#2 \dashrightarrow #3}}

\title{Abelianity conjecture for special compact K\" ahler threefolds}
\date{\today}
\author{Fr\'ederic Campana, Beno\^it Claudon}
\address{Fr\'ederic \textsc{Campana}, Beno\^it \textsc{Claudon}, Institut \'Elie Cartan Nancy, Universit\'e Henri Poincar\'e Nancy 1, B.P. 70239, 54506 Vandoeuvre-l\`es-Nancy Cedex, France}
\email{Frederic.Campana@iecn.u-nancy.fr, Benoit.Claudon@iecn.u-nancy.fr}

\begin{document}

\maketitle

\begin{flushright}
{\it Dedicated to V. Shokurov.}
\end{flushright}

\begin{abstract}
Using orbifold metrics of the appropriately signed Ricci curvature on orbifolds with negative or numerically trivial canonical bundle and the two-dimensional Log Minimal Model Program, we prove that the fundamental group of special compact K\"ahler threefolds is almost abelian. This property was conjectured in all dimensions in \cite{Ca04}, and also for orbifolds in \cite{Ca07}, where the notion of specialness was introduced. We briefly recall below the definition, basic properties, and the role of special manifolds in classification theory.
\end{abstract}

\section{Introduction}

We denote here by $X$ an $n$-dimensional compact connected K\"ahler manifold. The motivations for the following shortest, but non-transparent, definition of specialness will be explained below. See \cite{Ca04} where this notion was introduced for more details.

\begin{defi} A compact Kähler manifold $X$ is said to be \textbf{special} if, for any $p>0$, any rank-one coherent subsheaf $\lie\subset \Omega^p_X$ and any positive integer $N$, the natural meromorphic map $\merom{\Phi_{N.\lie}}{X}{\PP(V_{N,\lie}^*)}$ has an image of dimension at most $(p-1)$\footnote{By convention, this dimension is $-\infty$ if $V_{N,\lie}=\{0\}$ and is always at most equal to $p$ by a classical result of F. Bogomolov.}. Here $V_{N,\lie}^*$ denotes the dual of the complex vector space of sections of $Sym^N(\Omega^p_X)$ which take values in $\lie^{\otimes N}\subset Sym^N(\Omega^p_X)$ at the generic point of $X$.
\end{defi}

In particular, $X$ has then no surjective meromorphic map $\merom{f}{X}{Y}$ onto a manifold $Y$ of general type and dimension $p>0$, since otherwise $\lie:=f^*(K_Y)\subset \Omega^p_X$ would contradict the bound $(p-1)$ above. On the non-algebraic side, if $X$ has no map onto a positive dimensional-projective manifold, it is obviously special too. An equivalent  geometric definition of specialness of $X$ actually requires that $X$ satisfies the more restrictive condition of having no meromorphic map onto an orbifold of general type (see \ref{var speciale} for the precise definition and the relevant notions concerning orbifolds).

Specialness is preserved by bimeromorphic maps and finite \'elale covers, this last assertion being surprisingly difficult to show. Special manifolds generalise in higher dimensions the rational and elliptic curves,  which are obviously exactly the special curves. The next fundamental examples of special manifolds are indeed those which are either rationally connected or whith zero canonical dimension (\emph{i.e.} with $\kappa=0$, see example \ref{ex varietes speciales} below). However the class of special manifolds is much larger than the union of these two classes, since one shows by classification that a compact K\"ahler surface $X$ is special if and only if $\kod X\leq 1$ and $\gf X$ is almost abelian. In particular, ruled elliptic surfaces are special and surfaces with $\kod X=1$ are special if and only if they do not map onto any hyperbolic curve after some finite \'etale cover. No such simple characterisation is true when $n\geq 3$.

The central role of special manifolds in classification theory comes from the fact that, as shown in \cite[5.8]{Ca04}, any compact Kähler manifold is canonically and functorially decomposed by its \emph{core fibration} $c_X:X\to C(X)$ into its special part (the fibres of $c_X$) and its general type part (the orbifold base $(C(X),\Delta(c_X))$, which is its usual base $C(X)$ together with a ramification divisor $\Delta(c_X)$ on $C(X)$ encoding the multiple fibres of $c_X$).

Special manifolds and general type orbifolds are thus the {\it two} antithetical building blocks from which arbitrary compact K\"ahler manifolds are built in {\it one single step}. In contrast to  general type manifolds, for which no classification scheme  seems to be known or even expected, special manifolds are conjectured to have many fundamental properties in common with rational and elliptic curves.

Conjecturally indeed, an orbifold version of the $C_{n,m}$ conjecture implies that any special manifold is canonically and functorially decomposed, by means of orbifold versions of the rational quotient and of the Iitaka-Moishezon fibration, as a tower of fibrations whose orbifold fibres have either $\kappa=0$, or $\kappa_+=-\infty$ (a weak version of rational-connectedness, see \cite{Ca07}). Let us stress that the orbifold considerations are essential here (as in the LMMP, but for different reasons), and that apparently there is no possibility to work in the category of varieties without additional structure. This tower decomposition permits to lift (conditionally) to special manifolds\footnote{And even, more naturally, to special orbifolds.} properties which are expected to be common to manifolds which are either rationally connected or with $\kappa=0$ and naturally leads to the following conjectures:

\begin{conj}[\cite{Ca04,Ca07}]\label{conj abelianite}
1. Abelianity Conjecture: A special compact Kähler manifold has an almost abelian fundamental group (i.e: $\gf X$ has an abelian subgroup of finite index).\\
\noindent 2. A compact Kähler manifold (resp. a projective manifold defined over a number field) has an identically vanishing Kobayashi pseudometric (resp. is potentially dense) if and only if it is special (this last statement is inspired by Lang's conjectures).
\end{conj}

\noindent For example, rationally connected manifold are simply connected \cite{Ca92} and compact Kähler manifolds $X$ with $c_1(X)=0$ have an almost abelian fundamental group by \cite{Yau78} and \cite{Be83}.

In this article we prove the Abelianity conjecture above for threefolds (as said above, this is known for surfaces by classification, see Proposition \ref{surface speciale}).

\begin{thm}\label{main}
Let $X$ be a compact Kähler threefold. If $X$ is special, its fundamental group is almost abelian.
\end{thm}

This immediately implies, among several other things, a precise solution of Shafarevich's conjecture in this case (and for all special manifolds if the Abelianity Conjecture holds):

\begin{cor} If $X$ is a special compact K\" ahler threefold with universal cover $\wtx$, let $X'$ be any finite \'etale cover of $X$ with abelian torsionfree fundamental group and Albanese variety $Alb(X')$ of dimension $d:=q(X')$.Then $\wtx =X'\times_{Alb(X')}\widetilde{Alb(X')}$ is holomorphically convex, the universal cover $\widetilde{Alb(X')}$ of $Alb(X')$ being Stein, since isomorphic to  $\CC^d$.
\end{cor}

\section{Reduction to the 2-dimensional orbifold case}
We shall prove the theorem in this section only in the cases where no orbifold structure is needed, that is except when either $X$ is projective and $\kod X=2$, or when $a(X)=2$. The treatment of these two residual cases needs the consideration of two-dimensional projective special orbifolds and will be the subject of the subsequent sections.

Because of the lack of a Minimal Model program in the K\" ahler non-projective case, we need to treat it differently from the projective case. So assume first that $X$ is projective. We work according to the value of $\kod X\leq 2$.\\

\mathversion{bold}$\kod X=-\infty:$\mathversion{normal} By Miyaoka's theorem, $X$ is uniruled. Let $r_X:X\to R(X)$ be its rational quotient (also known as its MRC fibration). Then $R(X)$ is special too, with $\mathrm{dim}(R(X))\leq 2$, and $\gf X\simeq \gf{R(X)}$, since the fibres of $r_X$ are rationally connected. Since $\gf {R(X)}$ is almost abelian, so is $\gf X$.

\mathversion{bold}$\kod X=0:$\mathversion{normal} If $c_1(X)=0$, the theorem is true, by \cite{Yau78} and \cite{Be83}. One reduces to this case by the Minimal Model Program and \cite{NS95} (see \cite[(4.17.3)]{K95}). The details are as follows: there exists a terminal model $X'$ birational to $X$ such that $K_{X'}$ is torsion, hence trivial after finite \'etale (in codimension $1)$ cover. Then $X'$ has only $cDV$ singularities. It is thus smoothable in the projective category by \cite{NS95}. The conclusion follows, since $\gf X$ drops by specialisation.

\mathversion{bold}$\kod X=1:$\mathversion{normal} Let $J_X:X\to Y$ be the Moishezon-Iitaka fibration. Because its base is a curve, there exists a finite \'etale cover of $X$ (still written $X)$ such that $J_X$ has no multiple fibre\footnote{Or $Y\cong \PP^1$ and $J_X$ has one or two multiple fibres. This case is easily treated similarly.}. Then $\gf X$ is an extension of $\gf Y$ $(Y$ being a rational or elliptic curve, since special) by a quotient of $\gf{X_y}$, $X_y$ being the generic smooth fibre of $J_X$, which is a surface with $\kappa=0$. Because these two groups are almost abelian, $\gf X$ is {\it a priori} only polycyclic. The following result however implies the conclusion.

\begin{thm}[\cite{Ca01}]\label{tfbspec'}
Let $f:X\to Y$ be a fibration without multiple fibres in codimension $1$ on $Y$ from a compact K\" ahler manifold onto a manifold $Y$. Assume that $Y$ and the generic fibre $X_y$ of $f$ both have an almost abelian fundamental group. Then $X$ also has an almost abelian fundamental group.
\end{thm}

The proof of this theorem rests on two deep results of Hodge theory: Deligne's strictness theorem for morphisms of MHS, and Hain's functorial MHS on the Mal\v cev completion of $\gf X$ when $X$ is compact K\" ahler.

\mathversion{bold}$\kod X=2:$\mathversion{normal} When $X$ is projective, this is the only remaining case. Notice that when $\kod X=1$, we could remove the multiple fibres of $J_X$ by making a suitable finite \'etale cover of $X$ because $Y$ was a curve. When $\kod X=2$ this is no longer possible in general. This is the reason why the notion of orbifold base, which virtually removes the multiple fibres, is introduced in this case, and why the geometry of such orbifolds needs to be developed. But once this is done, and the corresponding properties established, the proof is entirely parallel.\\

We now deal with the case when $X$ is {\it not} projective. We work this time according to the algebraic dimension $a(X)\in \{0,1,2\}$ of $X$.\\

\mathversion{bold}$a(X)=0:$\mathversion{normal} The irregularity of $X$ can take the values $q(X)=0,\,2$ or $3$ (since $q(X)=1$ would imply $a(X)\geq 1)$. The assumption $a(X)=q(X)=0$ leads to the finiteness of $\pi_1(X)$ according to \cite[Cor. 5.7]{Ca95}. Recall that the Albanese map is connected and surjective when $a(X)=0$. Thus, If $q(X)=3$, then $X$ is birational to its Albanese variety, and $\gf X$ is Abelian. When $q(X)=2$, the Albanese fibration $\alpha_X:X\To \mathrm{Alb}(X)$ has no multiple fibre in codimension one by \cite[Prop. 5.3]{Ca04}. Moreover, the general fibre of $\alpha_X$ is elliptic or rational, by \cite[Th. 13.8]{Ue75}. From Theorem \ref{tfbspec'}, we conclude that $\pi_1(X)$ is almost abelian.

\mathversion{bold}$a(X)=1:$\mathversion{normal} The algebraic reduction $a_X:X\to A(X)$ is then a fibration onto a curve $A(X)$ with general fibre special \cite[Th. 2.39]{Ca04}; the fundamental groups of the base and of the fibre are then almost abelian (see Theorem \ref{surface speciale}). As in the projective case with $\kod X=1$, we may assume that there are no multiple fibres. Then a final application of Theorem \ref{tfbspec'} implies the result.\\

We are thus left with the following two cases:
\begin{enumerate}[(1)]
\item $X$ is projective and $\kod X=2$,
\item $X$ is a compact K\"ahler threefold with $a(X)=2$.
\end{enumerate}
Let us briefly explain how the conclusion is obtained there. See next sections for the notions which are being used. In both cases we have a fibration (which may be assumed to be neat after suitable modifications of $X$ and $S$) $f:X\to S$ on a smooth projective surface $S$, with generic fibres elliptic curves. This defines a smooth (hence klt) orbifold base $(S,\Delta_f)$. Because $X$ is special, so is $(S,\Delta_f)$. As above, for such a fibration, $\gf X$ is now an extension of $\gf{S,\Delta_f}$ by $\gf{X_s}$. We show below:
\begin{thm}\label{groupe orbisurface spéciale}
The fundamental group of a special geometric orbifold $(S,\Delta_S)$ of dimension 2 is almost abelian of even rank at most $4$.
\end{thm}
\noindent The conclusion now follows from the orbifold version of theorem \ref{tfbspec'} above.\\

Let us stress that our proof of theorem \ref{groupe orbisurface spéciale} depends in an essential way on the fact that \emph{integral klt} pairs are locally uniformisable by smooth germs of surfaces, because we are using the existence of Ricci-flat orbifold metrics when $c_1(S,\Delta)=0$. This is the main reason why we cannot extend the above theorem \ref{groupe orbisurface spéciale} to any higher dimension. Otherwise, using the Abundance Conjecture, it seems that one could derive the Abelianity Conjecture in the projective case for klt orbifolds in any dimension using inductively on the dimension essentially the same arguments as the ones presented below.

\section{Some basic facts on orbifolds}

\subsection{Notion of orbifold}
The orbifolds considered here are spaces with local smooth uniformisations under the action of some finite groups. To deal with fixed point in codimension one, the structure is enriched with a $\QQ$-divisor of a specific form. The definition below is taken from \cite{GK07}.
\begin{defi}\label{defi-orbi}
An orbifold is a pair $(X,\Delta)$ where $X$ is a normal variety and $\Delta$ a $\QQ$-Weil divisor of the following form:
$$\Delta=\sum_{i\in I}(1-\frac{1}{m_i})\Delta_i,$$
where the $m_i\ge2$ are integers, and the pair $(X,\Delta)$ is locally uniformizable in the following sense: $X$ is covered by the domains $U$ of finite maps
$$\varphi:U\To X$$
satisfying the following properties:
\begin{enumerate}[(i)]
\item $\varphi(U)$ is open in $X$,
\item $\varphi:U\To \varphi(U)$ is a Galois cover whose branching divisor is exactly $\Delta_{\vert\varphi(U)}$.
\end{enumerate}
The support of $\Delta$ is then: $\abs{\Delta}=\sum_i \Delta_i$.
\end{defi}

\noindent{\bf Terminological remark.}\\
The orbifolds we consider are compatible with all situations in which this term is used: they are particular cases of the ones in \cite{Ca07} and the special manifolds defined below are actually the \emph{classical} special ones in \emph{ibid.} They are smooth Deligne-Mumford stacks too and also klt pairs of the LMMP. Because they are locally smoothly uniformised, we can attach to them fundamental groups and differential-geometric notions such as metrics and differential forms.

\begin{defi}\label{fibré canonique orbi}
The canonical divisor of such a pair is the (Weil) $\QQ$-divisor $K_X+\Delta$.
\end{defi}

\begin{ex}\label{exemple snc}
Let $X$ be a smooth variety, and  $\abs{\Delta}=\sum_i \Delta_i$ be a normal crossing divisor; the choice of multiplicities $m_i\geq 2$ on each component of the divisor defines a canonical orbifold structure on $(X,\Delta)$. Since $\abs{\Delta}$ is locally given by the equation $z_1\cdots z_k=0$ in suitable coordinates $(z_1,\dots,z_n)$, the map $(z_1,\dots,z_n)\mapsto (z_1^{m_1},\dots,z_k^{m_k},z_{k+1},\dots,z_n)$
gives a local uniformization. These orbifolds are said to be smooth and integral in \cite{Ca07}.

In particular, an orbifold curve is simply a smooth curve with a finite set of points marked with integral multiplicities at least $2$.
\end{ex}
\noindent Since we consider only integral and finite multiplicities here, we define:
\begin{defi}\label{defi orbi geo}
A geometric orbifold is an orbifold $(X,\Delta)$ with $X$ smooth and $\abs{\Delta}$ of normal crossings.
\end{defi}

\subsection{Orbifold base of a fibration and special manifolds}

\begin{defi}[Def. 1.2, \cite{Ca04}]\label{defi fibration nette}
Let $f_0:X_0\To Y_0$ be a fibration (surjective morphism with connected fibres) between compact complex manifolds. A \emph{neat} model of $f$ consists in a commutative diagram:
$$\xymatrix{X\ar[d]_f \ar[r]^u & X_0\ar[d]^{f_0} \\
Y\ar[r]_v & Y_0
}$$
where
\begin{enumerate}[1.]
\item $X$ and $Y$ are smooth,
\item $u$ and $v$ are bimeromorphic morphisms,
\item the locus of singular fibres of $f$ is a normal crossing divisor of $Y$,
\item every $f$-exceptional divisor is $u$-exceptional too.
\end{enumerate}
\end{defi}
Such models actually exist birationally for any fibration; this can be proven using Hironaka desingularization and Raynaud flattening's theorems \cite[Lemma 1.3]{Ca04}. Notice that any fibration is neat when $Y$ is a curve.\\

Given a neat fibration $f:X\To Y$, we can naturally associate a $\QQ$-divisor $\Delta^*(f)$ on the base $Y$ of $f$. This divisor will be supported on the singular locus of $f$ and will thus be normal crossings: the pair $(Y,\Delta^*(f))$ will be a geometric orbifold.

The construction goes as follow (see also \cite[1.1.4]{Ca04}): if $\Delta_i$ is any component of the singular locus of $f$, its pull-back can be written
$$f^*(\Delta_i)=\sum_{j}m_j D_j+R$$
where $R$ is $f$-exceptional and $D_j$ is mapped surjectively onto $\Delta_i$. The multiplicity of $f$ along $\Delta_i$ is defined by:
$m_i=m(f,\Delta_i)=\mathrm{gcd}_j(m_j).$
\begin{defi}\label{base orbi stable}
The pair $(Y,\Delta^*(f))$, where:
$\Delta^*(f)=\sum_i (1-\frac{1}{m_i})\Delta_i,$
is called the \emph{orbifold base} of the fibration $f_0$ (notations as above).\\
A fibration $f_0$ is said to be of general type if the canonical divisor of the orbifold base $(Y,\Delta^*(f))$ is \emph{big}:
$\kappa(Y,K_Y+\Delta^*(f))=\dimm{Y}>0.$
\end{defi}

\begin{defi}[\cite{Ca04}, Def. 2.1]\label{var speciale}
A compact Kähler manifold $X$ is said to be (classically) special if it does not admit any fibration of general type.
\end{defi}

\begin{ex}\label{ex varietes speciales}
The main examples of special manifolds are given by the following classes \cite[3.22, 5.1, 2.39]{Ca04}:
\begin{enumerate}[$\star$]
\item rationally connected manifolds,
\item compact Kähler manifolds $X$ with $\kappa(X)=0$. This is a consequence of the additivity of canonical dimensions in general type fibrations (\cite[Th. 4.2]{Ca04}, an orbifold version of Viehweg's theorem).
\item The special curves are thus just the rational or elliptic ones.
\item compact K\"ahler manifolds of algebraic dimension zero and more generally fibres of algebraic reductions.
\end{enumerate}

\end{ex}
Conjecturally in an orbifold version of Iitaka's $C_{n,m}$-conjecture, special manifolds can be reconstructed as a tower (in a suitable sense) of fibrations with fibres belonging to the classes above. This reduces Abelianity Conjecture \ref{conj abelianite} to the case of orbifolds with either $\kappa=0$ or $\kappa_+=-\infty$ (see \cite[13.10]{Ca07}). 

In dimension $2$ only, we still have a simple topological characterisation of specialness:

\begin{prop}[Prop. 3.32 ,\cite{Ca04}]\label{surface speciale}
A compact Kähler surface $X$ is special if and only if: $\kappa(X)\leq 1$ and if $\pi_1(X)$ is almost abelian.
\end{prop}

\noindent We shall need also the notion of orbifold base of a fibration $f:X\to Y$ also when $X$ is equipped with an orbifold divisor $\Delta_X$, at least when $Y$ is a curve and $X$ is a surface (for the general case, see the definition in \cite[Def. 4.2]{Ca07}).
\begin{defi}
Let $(X,\Delta_X)$ be a geometric orbifold of dimension 2 and $f:X\To C$ a fibration onto a curve. Let us define the multiplicity of a point $y\in C$ (relatively to $f$ and $\Delta_X$) by the following formula:
$m_y(f,\Delta_X):=\mathrm{gcd}_i\left(m_i. \mathrm{mult}_{\Delta_X}(F_i)\right),$
where $f^*(y)=\sum_i m_iF_i$.
The orbifold base is then the pair $(C,\Delta^*(f,\Delta_X))$ where:
$$\Delta^*(f,\Delta_X)=\sum_y (1-\frac{1}{m_y(f,\Delta_X)}).\set{y}$$
We say that  $f$ is of general type if:
$2g(C)-2+\mathrm{deg}\left(\Delta^*(f,\Delta_X)\right)>0.$

\noindent The $2$-dimensional orbifold $(X,\Delta_X)$ is then said to be special if $(X,\Delta_X)$ does not admit any general type fibration onto a curve and if $\kappa(X,K_X+\Delta_X)<2$.
\end{defi}

\section{Fundamental groups and fibrations}

\subsection{An orbifold exact sequence}

In this paragraph, we define fundamental groups of orbifolds and study the morphisms induced at the level of fundamental groups by classical orbifold morphisms. 

\begin{defi}\label{groupe fondamental}
Let $(X,\Delta)$ be an orbifold; its fundamental group $\pi_1(X,\Delta)$ is the quotient of the group $\pi_1\left(X^*\backslash \abs{\Delta}\right)$ by the normal subgroup generated by the loops $\gamma_j^{m_j}$ where $\gamma_j$ is a small loop around the component $\Delta_j$ of multiplicity $m_j$ (where $X^*$ denote the smooth locus of $X$).
\end{defi}
This definition is derived from the local models. If $X=\CC^n/G$ where $G$ is a finite subgroup of $\mathrm{GL}_n(\CC)$ and $\Delta$ is the branching divisor of the projection $\pi:\CC^n\To X$, we recover the group $G$ as an orbifold fundamental group\footnote{The orbifold divisor determines the structure of the smooth Deligne-Mumford stack $\XX$ associated with $(X,\Delta)$ and called the \emph{root stack} by Abramovich-Vistoli.}:
\begin{equation*}\label{groupe fond local}
\pi_1(X,\Delta)\simeq G.
\end{equation*}

\begin{ex}\label{exemple groupe courbe}
The case of orbifold curves is quite classical. The structure of the fundamental group of such an orbifold curve $(C,\Delta=(m_1,\dots,m_n))$ (we just keep in mind the deformation invariant, that is: the multiplicities, not the marked points) is determined by the sign of its canonical bundle:
$$\mathrm{deg}\left(K_C+\Delta\right)=2g(C)-2+\sum_{j=1}^n (1-\frac{1}{m_j}).$$
When this quantity is positive (resp. zero, resp. negative), the fundamental group of $(C,\Delta)$ is commensurable to the fundamental group of an hyperbolic curve (resp. commensurable to $\ZZ^2$, resp. finite and explicitely known).
\end{ex}

When dealing with fibrations having multiple fibres, we need to consider orbifold fundamental groups. A fibration $f:X\To Y$ with general fibre $X_y$ gives rise to  a natural sequence of fundamental groups: 
$$\pi_1(X_y)\stackrel{i_*}{\To} \pi_1(X)\stackrel{f_*}{\To}  \pi_1(Y)\To 1.$$
Although $f_*$ is surjective since $f$ has connected fibres, this sequence is in general not exact in the middle. Multiple fibres are responsible for this failure, remedied by the orbifold fundamental group.

We will need the case where $X$ has an orbifold structure and so a slightly more general version. For the definition (of a neat fibration in the orbifold setting), we refer to \cite[Def. 4.8]{Ca07}. 

\begin{prop}[\cite{Ca07}, Cor. 12.10]\label{suite exacte}
Let $f:(X,\Delta_X)\To Y$ be a neat fibration. If $X_y$ denote a general fibre of $f$ and $\Delta_y$ the restriction of the orbifold structure to $X_y$, the sequence
$$\pi_1(X_y,\Delta_y)\To \pi_1(X,\Delta_X)\stackrel{f_*}{\To}  \pi_1(Y,\Delta^*(f,\Delta_X))\To 1$$
is exact.
\end{prop}

\begin{rem}\label{remarque fibration nette}
We omit the detailled definition of neatness required for the preceding statement because in the sequel we will only use it in the two following quite simple situations:
\begin{enumerate}[(a)]
\item when $\Delta_X$ is empty, this neatness assumption has been already encountered (\emph{cf.} Definition \ref{defi fibration nette}). The content of the previous proposition is then that $\pi_1(X)$ sits then in the middle of a short exact sequence (with the fundamental group of the orbifold base on the right hand).
\item when $Y$ is a curve: the fibration $f:(X,\Delta_X)\To Y$ is always neat.
\end{enumerate}
\end{rem}

\subsection{Nilpotency class in fibrations}

To complete this study of the behaviour of fundamental groups in fibrations, we prove that nilpotency conditions on the fundamental groups are preserved in fibrations between Kähler orbifolds. This remarkable fact (obviously false even for submersions between complex manifolds as the Iwasawa manifold shows) is a consequence of deep results in Mixed Hodge Theory\footnote{The key ingredients are: existence of a \textsc{mhs} on the Mal\v{c}ev completion of the fundamental group of a compact Kähler manifold and the strictness of morphisms of \textsc{mhs}.}.

\begin{thm}[Cor. 7.6, \cite{Ca10}]\label{fibration presque abélien}
Let $f:(X,\Delta_X)\To (Y,\Delta_Y)$ be a neat fibration (see Remark \ref{remarque fibration nette}) between smooth compact Kähler orbifolds. If the groups $\pi_1(X_y,\Delta_{X_y})$ and $\pi_1(Y,\Delta_Y)$ are almost abelian\footnote{More generally, if these two groups are torsionfree nilpotent of nilpotency class at most $\nu$, then $\pi_1(X,\Delta_X)$ is also nilpotent of nilpotency class at most $\nu$.} , then $\pi_1(X,\Delta_X)$ is almost abelian as well.
\end{thm}
\begin{demo}
It is a reduction to the case when $\Delta_X=0$. Since $f$ is assumed to be neat, the sequence
$$\pi_1(X_y,\Delta_{X_y})\To \pi_1(X,\Delta_X)\stackrel{f_*}{\To}  \pi_1(Y,\Delta_Y)\To 1$$
is exact. Since $\pi_1(X_y,\Delta_{X_y})$ and $\pi_1(Y,\Delta_Y)$ are almost abelian groups (of finite type), $G=\pi_1(X,\Delta_X)$ is then an almost polycyclic group; in particular, the group $G$ is linear. A famous result of Selberg asserts that $G$ has then a finite index subgroup $G'\le G$ which is torsion free. Consider $X'$ the orbifold cover of $(X,\Delta_X)$ associated with $G'$; the latter being torsion free, $X'$ is a normal variety with no orbifold structure (\emph{i.e.} $\Delta'=0$). Since $X$ has only quotient singularities, its fundamental group is isomorphic to the one of $\tilde{X}$, a desingularisation of $X'$ \cite[Th. 7.5]{K93}. To conclude, we consider (a neat model of) the Stein factorization of
$$\tilde{X}\To X'\To X\To Y.$$
Indeed, taking further blow-up of $\tilde{X}$, we can complete the picture:
$$\xymatrix{\tilde{X}\ar[r]\ar[rrd]_{\tilde{f}} & X'\ar[r] & X\ar[r]^{f} & Y\\
&&\tilde{Y}\ar[ru]_{g}&
}$$
with $\tilde{f}$ a neat fibration and $g$ generically finite. It is then easy to see that we have an exact subsequence (with vertical maps having finite index images):
$$\xymatrix{\pi_1(\tilde{X}_y) \ar[r] \ar@{^{(}->}[d] & \pi_1(\tilde{X})  \ar[r]^{\tilde{f}_*} \ar@{^{(}->}[d] &  \pi_1(\tilde{Y},\Delta^*(\tilde{f}))\ar[r] \ar@{^{(}->}[d] &1\\
\pi_1(X_y,\Delta_{X_y}) \ar[r] & \pi_1(X,\Delta_X)\ar[r]^{f_*} & \pi_1(Y,\Delta_Y)\ar[r] &1
}$$
So we are reduced to the same situation with $\Delta_X=0$ and we can apply \cite[Cor. 7.6]{Ca10}.
\end{demo}

\subsection{Fundamental groups of orbifolds with non-positive canonical bundle}

As for of compact Kähler manifolds, when the first Chern class is zero or positive, the now usual differential-geometric methods can then be applied to construct orbifold K\"ahler metrics with Ricci-curvature of the corresponding sign (using the local smooth uniformisations).

\begin{thm}\label{groupe orbi c1 positif}
Let $\XX=(X,\Delta)$ be an orbifold with $X$ compact K\" ahler. If the first Chern class of $\XX$ is non negative, the group $\pi_1(\XX)$ is almost abelian; more precisely:
\begin{enumerate}[$\star$]
\item if $c_1(K_X+\Delta)=0$, then $\pi_1(X,\Delta)$ is almost abelian of even rank\footnote{The rank of an almost abelian group $G$ of finite type is the maximum rank of an abelian  subgroup of finite index of $G$.} bounded by $2\dimm{X}$.
\item if $c_1(K_X+\Delta)>0$, then $\pi_1(X,\Delta)$ is finite.
\end{enumerate}
\end{thm}

To begin with, let us recall some basic facts on differential calculus on orbifolds. The smooth functions (differential forms, hermitian metrics,$\dots$) on an orbifold $(X,\Delta)$ are the smooth functions on $X^*\backslash\abs{\Delta}$ which can be smoothly extended (as the usual objects, after taking inverse images) in local uniformizations. For instance, if $(X,\Delta)$ is a geometric orbifold, a Kähler metric has the following form in coordinates charts adapted to $\Delta$:
\begin{align*}
\omega_{\Delta}&=\omega_{eucl}+\sum_{j=1}^n i\partial\overline{\partial}\abs{z_j}^{2/m_j}\\
&=\omega_{eucl}+\sum_{j=1}^n\frac{idz_j\wedge d\overline{z_j}}{m_j^2\abs{z_j}^{2(1-1/m_j)}}.
\end{align*}
Local uniformizations can also be used to compute integrals of forms of maximal degree. If $\varphi:U\To \varphi(U)\subset X$ is such a local cover and $\alpha$ an orbifold top-form on $X$, define
$$\int_{\varphi(U)}\alpha=\frac{1}{\mathrm{deg}(\varphi)}\int_U \varphi^*\alpha.$$
This local computation can then easily be globalized using partitions of unity.

The canonical divisor has been defined in \ref{fibré canonique orbi} as a $\QQ$-divisor. It should be noted here that ($X$ being compact and the uniformizations being finite) some integral multiple of this divisor defines a line bundle on $X$ and it can be used to compute the first Chern class of $X$. As a byproduct, the Ricci form of any orbifold volume form is a (1,1)-orbifold form whose cohomology class coincides with $c_1(\XX)$.\\

As in the manifold case, every invariant form whose cohomology class coincides with $c_1(\XX)$ is the Ricci curvature of an orbifold K\"ahler metric. This fact has been already noticed several times in the literature (see for instance \cite{J00,Ca4} when $c_1(\XX)=0$, and \cite{DK01} when $c_1(\XX)>0$, and the references therein). For sake of completeness, let us recall the statement.
\begin{thm}[Calabi-Yau]\label{orbi-calabi}
Let \XX be a K\"ahler orbifold whose underlying space is compact and let us fix $\omega_0$ an orbifold K\"ahler metric. For any representative (\emph{i.e.} smooth invariant $(1,1)$ form) $\alpha$ of $c_1(\XX)$, there exists a unique orbifold K\"ahler metric $\omega$ in the K\"ahler class $[\omega_0]$ such that:  $\mathrm{Ricci}(\omega)=\alpha.$
\end{thm}
\begin{demo}
We reproduce the arguments given in \cite[Th. 4.1]{Ca4}. This is a simple adaptation of the proof exposed in \cite[p. 85-113]{S89}. Since $\mathrm{Ricci}(\omega_0)$ and $\alpha$ define the same class, the orbifold $\partial\overline{\partial}$-lemma provides us a smooth (orbifold) function $f$ such that $\mathrm{Ricci}(\omega_0)=\alpha+i\partial\overline{\partial}f$. We normalize $f$ (by addition of a suitable constant) in such a way that
$$\int_X (e^f-1)\omega_0^n=0.$$
The problem is then reduced to solving the Monge-Ampère equation
$$(\mathrm{MA})\quad\left(\omega_0+i\partial\overline{\partial}\varphi\right)^n=e^f\omega_0^n.$$
To do so, we applied the continuity method as in \cite[p. 85-113]{S89} and consider the set $T$ of $t\in [0,1]$ for which the equation
$$(\mathrm{MA})_t\quad\left(\omega_0+i\partial\overline{\partial}\varphi_t\right)^n=C_te^{tf}\omega_0^n.$$
has a solution where
$$C_t=\frac{\int_X\omega_0^n}{\int_Xe^{tf}\omega_0^n}.$$
When $t=0$, $\varphi_0=0$ is an obvious solution for $(\mathrm{MA})_0$. The set $T$ is open by linearization of the problem and implicit functions theorem (applied in local uniformizations). To see that $T$ is closed we need a priori estimates on the solution of $(\mathrm{MA})_t$ which are established using two types of arguments: maximum principle and Nash-Moser iteration (integral inequalities). These can be used in the orbifold setting as well: maximum principle is applied in local uniformizations and integral inequalities are also valid in this formalism (see above).
\end{demo}

We proceed as in the manifold case to deduce Theorem \ref{groupe orbi c1 positif} from Theorem \ref{orbi-calabi}. The positivity assumptions can be translated in the following way:  in the Fano case $(c_1(\XX)>0)$ there exists a positive $(1,1)$ form $\alpha$ representing $c_1(\XX)$;  in the Ricci flat case $(c_1(\XX)=0)$ we choose $\alpha=0$. We can now apply Theorem \ref{orbi-calabi}:  \XX has an orbifold K\"ahler metric $\omega$  with prescribed Ricci curvature: $\mathrm{Ricci}(\omega)=\alpha.$

In the Fano case, by compactness of $X$, we get: $\mathrm{Ricci}(\omega)\geq \epsilon\omega$,
where $\epsilon>0$. The following orbifold version of Myer's theorem applied to the orbifold universal cover\footnote{see \cite{Na87}.} of $\XX$ implies the finiteness of the orbifold fundamental group: 

\begin{thm}[\cite{Bor93}]\label{orbi-myer}
Let $(\mathcal{Y},g)$ be a complete Riemannian orbifold (of dimension $n$). If the Ricci curvature satisfies the inequality
$$\mathrm{Ricci}(g)\ge(n-1)k$$
for some $k>0$, the underlying space $Y$ of $\mathcal{Y}$ is then compact with diameter bounded above:
$$\mathrm{diam}(Y)\le\frac{\pi}{\sqrt{k}}.$$
\end{thm}

When $c_1(\XX)=0$, the following orbifold splitting theorem \ref{orbi decomposition} applied also to the orbifold universal cover of $\XX$ also implies the claim.

\begin{thm}[\cite{BZ94}]\label{orbi decomposition}
Let $(\mathcal{Y},g)$ be a compact Riemannian orbifold (of dimension $n$). If the Ricci curvature is everywhere non negative, the orbifold universal cover $\widetilde{\mathcal{Y}}$ admits a (metric) splitting
$$\widetilde{\mathcal{Y}}\simeq N\times\RR^m$$
where $N$ is a compact orbifold. The orbifold fundamental group of $\mathcal{Y}$ is also an extension of a crystallographic group by a finite group and is in particular almost abelian.
\end{thm}

From the arguments given in \cite[5.4, 6.3]{Ca4} we deduce that the rank of the almost group $\pi_1(X,\Delta)$ is even and bouded by $2\mathrm{dim}(X)$.$\square$

\section{Minimal Model Program for klt pairs in dimension $2$}

\subsection{klt pairs as orbifolds in dimension $2$}

In this short paragraph, we gather from \cite{KM98} some well known facts on the MMP for log pairs in dimension 2. Starting with a pair $(X,\Delta)$ where $X$ is a smooth surface and $\Delta=\sum_j b_j\Delta_j$ a $\QQ$-Weil effective divisor, we can perform a sequence of divisorial contractions which ends with a birational model of the initial pair and whose geometry is simplest, according to the sign of the canonical bundle. In this process, however, $(K_X+\Delta)$-negative curves are being contracted, and the resulting surface is no longer smooth in general. The relevant preserved category of singularities is then described as follows.
\begin{defi}\label{def klt}
A pair $(X,\Delta)$ where $X$ is a $\QQ$-factorial normal variety is said to have only klt singularities if
\begin{enumerate}[(i)]
\item $\forall\, j,\, 0<b_j<1$
\item For any (or equivalently, one) log-resolution $f:Y\To X$ of $(X,\Delta)$, in the following equality, we have $a_i>-1$, for all $i's$:
$$K_Y+\tilde{\Delta}=f^*(K_X+\Delta)+\sum_{i}a_i E_i,$$
$\tilde{\Delta}$ being the strict transform of $\Delta$, and $E_i$ the exceptional divisors of $f$.
\end{enumerate}
\end{defi}
A smooth (integral) orbifold is of course a $klt$ pair. These singularities are preserved in the LMMP:

\begin{thm}[Th. 3.47, \cite{KM98}]\label{LMMP surface}
Let $(X,\Delta)$ be a klt surface. There exists a birational morphism $f:X\To S$ such that the resulting pair $(S,\Delta_S=f_*(\Delta))$ is still klt and satisfies (exactly) one of the following properties :
\begin{enumerate}
\item $K_S+\Delta_S$ is nef,
\item $S$ admits a fibration $\pi:S\To C$ onto a (smooth) curve $C$, the general fibre of $\pi$ being a smooth $K_S+\Delta_S$-negative rational curve,
\item $(S,\Delta_S)$ is Del Pezzo : $\rho(S)=1$ and $-(K_S+\Delta_S)$ is ample.
\end{enumerate}
\end{thm}

When $\Delta=0$, it is well know that klt singularities coincide with quotient singularities in dimension 2 \cite[Prop. 4.18]{KM98}. This property still holds for $\Delta\neq 0$ with integral multiplicities (\emph{i.e.} coefficients $b_j$ of the form $1-\frac{1}{m}$).
\begin{thm}\label{klt2=quotient}
Let $(X,\Delta)$ be a pair where $X$ is a surface and $\Delta$ has integral multiplicities. The following conditions are equivalent near each point of $X$:

\begin{enumerate}
\item $(X,\Delta)$ is klt,
\item $(X,\Delta)$ has a finite local fundamental group,
\item $(X,\Delta)$ is locally presented as a quotient: $\pi:\CC^2\To \CC^2/G\simeq X$,
where $G$ is a finite group acting linearly on $\CC^2$ and $\Delta$ is the unique $\QQ$-divisor on $X$ such that $\pi^*(K_X+\Delta)=K_{\CC^2}$.
\end{enumerate}
\end{thm}
This statement is found and used in several places \cite{K90,M00}, but never with a complete accessible proof: the appendix of \cite{K90} consists in a list of the possible cases, and refers to the thesis \cite{Nak} for the proof. For this reason, we give a proof of Theorem \ref{klt2=quotient} at the end of the present text (\emph{cf.} Appendix).

\subsection{Fundamental groups and Mori contractions}

We now relate the fundamental groups of a geometric orbifold $(X,\Delta)$ of dimension 2 and of its minimal model $(S,D=\Delta_S=f_*(\Delta))$ as in theorem \ref{LMMP surface}.

\begin{prop}\label{groupe modele minimal}
Let $(X,\Delta)$ be a geometric orbifold of dimension 2 and $f:(X,\Delta)\To (S,D)$ its minimal model. There is a natural surjective morphism of groups:
$$f^\sharp:\pi_1(S,D)\To \pi_1(X,\Delta).$$
In particular, if $\pi_1(S,D)$ is almost abelian so is $\pi_1(X,\Delta)$.
\end{prop}
\begin{demo}
Let us call $E$ the union of the curves contracted by $f$; $f$ being an isomorphism away from $E$ we have natural maps, $S^*$ being the smooth locus of $S$:
$$S^*\backslash\abs{D}\stackrel{f^{-1}}{\To} X\backslash\left(\abs{\Delta}\cup E\right)\hookrightarrow X\backslash\abs{\Delta}$$
(note that $\abs{\Delta}$ and $E$ may have common components). At the level of fundamental groups, we get natural morphisms:
$$\pi_1(S^*\backslash\abs{D})\stackrel{\sim}{\To}\pi_1\left(X\backslash(\abs{\Delta}\cup E)\right)\twoheadrightarrow
\pi_1(X\backslash\abs{\Delta})\twoheadrightarrow\pi_1(X,\Delta).$$
To get the morphism $f^\sharp$ we need only to remark that the loops around the components of $D$ are sent onto loops around corresponding components in $\Delta$ with the same multiplicities, since $f$ is an isomorphism between $S^*$ and $X\backslash E$.
\end{demo}

\begin{ex}\label{attention} In general, the kernel of this morphism can however be very big as shown by the following example.

Let $C$ be an elliptic curve, and let $X$ be the blow-up of $S=C\times C$ at a point $(c,c)$. Let $F\subset S$ the fibre of the first projection through this point and let $G$ be its strict transform in $X$. If $m>1$ is an integer, the minimal model of $(X,\Delta=(1-\frac{1}{m})G)$ is $(S,D=(1-\frac{1}{m})F)$ and the exceptional divisor $E$ of the blow-up is the only $(K_X+\Delta)$-negative curve on $X$. Since the orbifold base of the first projection is $(C,(1-\frac{1}{m})\set{c})$, the fundamental group of $(S,D)$ has (a finite index subgroup which have) a surjective morphism onto a non-abelian free group. On the other hand the orbifold base of $g:(X,\Delta)\to C$, $g$ being the composition of the first projection and the blow-up, is merely $C$ since the fibre over $c$ has $E$ as component and inherits from it the multiplicity 1. It is then easy to see that $(X,\Delta)$ is special and that its fundamental group is isomorphic to the fundamental group of $X$ and thus abelian.
\end{ex}

\begin{rem}\label{remarque groupe fondamental mmp}
The direction of the arrow $f^\sharp$ in proposition \ref{groupe modele minimal} does not look to be functorial. So we explain its construction differently. The map $f:(X,\Delta)\to (S,D)=(X,\Delta)_{min}$ is in general not a `classical' orbifold morphism (in the sense of \cite{Ca07}) and thus does not induce a functorial morphism of groups. The map $f$ induces an orbifold morphism, and thus a morphism of groups, only when the multiplicities on the exceptional divisors of $f$ are sufficiently divisible (for example by the order of the local fundamental group of $(S,D)$ at the point under consideration). In our example \ref{attention} above, equipping $E$ with a multiplicity divisible by $m$ is the right condition. 

Let thus $\Delta^+$ be an orbifold divisor on $X$ such that $f_*(\Delta^+)=f_*(\Delta)=D$, and such that $f: (X,\Delta^+)\to (S,\Delta)$ and $id_X:(X,\Delta^+)\to (X,\Delta)$ are orbifold morphisms (this means for $id_X$ that the multiplicity of any component of $\Delta$ divides the corresponding multiplicity in $\Delta^+)$. 

We then get two functorial group morphisms: $f_+^\sharp:\pi_1(X,\Delta^+)\to \pi_1(S,D)$, which is an isomorphism, and $id_X^\sharp:\pi_1(X,\Delta^+)\to \pi_1(X,\Delta)$, which is surjective (as above). Our initial $f^\sharp$ was nothing but $(id_X)^\sharp\circ (f_+^\sharp)^{-1}$.
\end{rem}

\subsection{Abelianity for special klt pairs and proof of the main theorem}

To conclude this section, we shall use the preceding proposition to prove Abelianity Conjecture for geometric orbifolds of dimension 2.\\

\begin{demo}[of Theorem \ref{groupe orbisurface spéciale}]
Let $(S,D)$ the minimal model of the pair $(X,\Delta_X)$. We argue according to the values of $\kappa=\kod{X,\Delta}=\kod{S,D}\in\set{-\infty,0,1}$.

If $\kappa=1$, then $X$ admits a fibration $f$ onto a smooth curve $C$ whose general fiber $F$ satisfies $F\cdot (K_X+\Delta)=0$; the general orbifold fibre $(F,\Delta_F)$ is then orbifold-elliptic and its fundamental group is almost abelian (example \ref{exemple groupe courbe}). Adding to $C$ the orbifold divisor $\Delta^*=\Delta^*(f,\Delta_X)$ of classical multiplicities, we get an exact sequence
$$\pi_1(F,\Delta_F)\To \pi_1(X,\Delta_X)\To \pi_1(C,\Delta^*)\To 1.$$
The orbifold curve $(C,\Delta^*)$ is special since so is $(X,\Delta)$ and Theorem \ref{fibration presque abélien} shows that $\pi_1(X,\Delta)$ is almost abelian.

When $(S,D)$ admits a structure of Mori fibre space over a (special) curve (in particular $\kappa=-\infty$), we still have a fibration on $X$ whose fibres are special (in fact rational) and we can proceed as in the lines above.

There are two cases left:
\begin{enumerate}
\item $\kod{S,D}=0$,
\item or $(S,D)$ is log Del Pezzo.
\end{enumerate}
In the first case, $K_S+D$ being \emph{nef} it is semi-ample (log abundance for surfaces, see \cite{FMK}) and then torsion. In particular, the orbifold $(S,D)$ is Ricci-flat. In the second case, $(S,D)$ is Fano. In both cases, we can apply Theorem \ref{groupe orbi c1 positif} to conclude that the fundamental group of $(S,D)$ is almost abelian. Finally, Proposition \ref{groupe modele minimal} can be applied (it is only here we need it) to show that $\pi_1(X,\Delta_X)$ is almost abelian as well.
\end{demo}



\appendix

\section{Classification of klt singularities of surfaces}

In this appendix we give a proof of the classification of klt singularities for integral pairs (Theorem \ref{klt2=quotient}). The method is first to treat the singular case (using the method of \cite[Th. 4.7]{KM98}) and then to reduce the smooth case to the preceding one using orbifold \'etale covers. Let then $(X,\Delta)$, $\Delta=\sum_j(1-\frac{1}{m_j})D_j$ be a germ of klt pair. 

\subsection{The singular case.}

To begin with we quote a useful (negativity) lemma on connected quadratic forms.
\begin{lem}[Cor. 4.2, \cite{KM98}]\label{lfq}
Let $E:=\cup_j E_j$ be a connected exceptional curve on a smooth complex surface. Let $A:=\sum_j a_jE_j$ and $B:=\sum_j b_jE_j$, with $a_j,b_j\in \RR$. Assume that $A\cdot E_j\geq B\cdot E_j$ for any $j$. Then either $A=B$ or $a_j<b_j$ for any $j$.
\end{lem}

\noindent We first consider the case in which the germ $X$ is singular.

\begin{lem}
Assume that $X$ is singular, and that $\Delta\neq 0$. Let $f:X'\to X$ be the minimal resolution of the germ $X$, $\Delta'=\sum_k(1-\frac{1}{m_k})D'_k$ being the strict transform of $\Delta$ in $X'$ and $E=\cup_j E_j$ the exceptional divisor. The extended dual graph of $f^*(\Delta)$ is then one of the following ones below, in which the (non-compact) components of $\Delta'$ are indicated by black dots. Moreover, all intersections are transversal and all white dots are smooth rational curves. These are $(-2)$-curves, except possibly in the first case (see remark \ref{rklt+} below for additional constraints in the first case).

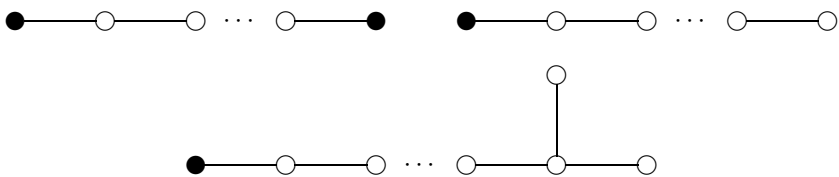
\begin{figure}[!h]
\setlength{\unitlength}{1.2mm}
\begin{picture}(120,18)
\put(5,17){\circle*{2}}
\multiput(15,17)(10,0){3}{\circle{2}}
\put(45,17){\circle*{2}}
\multiputlist(10,17)(10,0){{\line(1,0){8}},{\line(1,0){8}},{$\cdots$},{\line(1,0){8}}}
\put(55,17){\circle*{2}}
\multiput(65,17)(10,0){4}{\circle{2}}
\multiputlist(60,17)(10,0){{\line(1,0){8}},{\line(1,0){8}},{$\cdots$},{\line(1,0){8}}}

\put(25,1){\circle*{2}}
\multiput(35,1)(10,0){5}{\circle{2}}
\put(65,11){\circle{2}}\put(65,2){\line(0,1){8}}
\multiputlist(30,1)(10,0){{\line(1,0){8}},{\line(1,0){8}},{$\cdots$},{\line(1,0){8}},{\line(1,0){8}}}
\end{picture}
\caption{Dual graphs when $\Delta\neq0$}
\end{figure}

\end{lem}
\begin{demo} We assume the knowledge of the classification of Duval (\emph{i.e.} canonical) and of klt germs of $2$-dimensional singularities, and more precisely the fact that the dual graphs of their minimal resolutions are given either by Dynkin diagrams of $(-2)$-curves of type $A_n$, $D_n$, $E_k$ (for $k=6,7,8$) or by Hirzebruch-Jung chains with transversal intersections (\cite[Th. 4.7]{KM98}, the method of proof we shall now adapt).

From negativity, rational numbers $a_i$ are uniquely determined by the numerical equalities:
$$(\star)\quad\forall j,\,d_j:=(K_{X'}+\Delta')\cdot E_j=A\cdot E_j$$
with $A:=\sum_ia_iE_i$. Write also: $e_i:=-E_i^2$.


Assume first that there are (at least) two components $D_1$ and $D_2$ in $\Delta'$, of multiplicities $m_1,\,m_2$ . Consider a shortest chain $E_1,\dots,E_n$ of components of $E$ joining them, such that $D_1$ (resp. $D_2)$ meets $E_1$ (resp. $E_n)$, but no other $E_{\ell}$. An easy computation gives:
$$d_j=\left\{\begin{array}{l}(e_1-2)+(1-\frac{1}{m_1})(D_1\cdot E_1) \textrm{ if }j=1,\\
(e_n-2)+(1-\frac{1}{m_n})(D_2\cdot E_n) \textrm{ if }j=n,\\
e_j-2\textrm{ otherwise.}
\end{array}\right.$$
Consider now the curve $B:=-\beta(\sum_1^n E_{\ell})$, for some $\beta >0$ to be chosen latter. We assume first that $n\geq 2$. The intersection numbers are then given by:
$$B\cdot E_j=\left\{\begin{array}{l}\beta (e_j-1),\,j=1,n,\\
\beta (e_j-2),\, 1<j<n.
\end{array}\right.$$
We choose now the value of $\beta$ according to the self-intersection of the $E_j$:
\begin{enumerate}
\item if $e_i\neq 2$, for some $i\neq 1,n$, we set
\begin{align*}
\beta &=\inf \{1,\,\frac{d_j}{e_j-1},\,j=1,n\}\\
&=\inf \{1,\,1+\frac{(1-\frac{1}{m_1})(E_1\cdot D_1)-1}{e_1-1},\,1+\frac{(1-\frac{1}{m_2})(E_n\cdot D_2)-1}{e_n-1}\}
\end{align*}
\item and choose $\beta= \inf \{\frac{d_j}{e_j-1},j=1,n\}$ if $e_i=2,\forall i\neq 1,n$.
\end{enumerate}
We thus have $A\cdot E_j\geq B\cdot E_j$ for any $j$, and equality for some $j$. From Lemma \ref{lfq} we get $A=B$ and $a_j=-\beta>-1$ since $(X,\Delta)$ is klt. Using $(\star)$ it easily implies
$$(1-\frac{1}{m_1})(D_1\cdot E_1)<1\quad\textrm{and}\quad(1-\frac{1}{m_2})(D_2\cdot E_n)<1.$$
Because $(1-\frac{1}{m})\geq \frac{1}{2}$ we get $D_1\cdot E_1=D_2\cdot E_n=1$ and we are thus in the first case.
 
Assume now that $n=1$ and notice that then $D_i$ meets $E_1$ for $i=1,\,2$. We then have
$$(e_1-2)+(1-\frac{1}{m_1})(D_1\cdot E_1)+(1-\frac{1}{m_2})(D_2\cdot E_1)=\beta e_1^2$$ for some $\beta<1$ (by the klt condition). The solutions are easily determined
\begin{enumerate}[(i)]
\item either $E_1\cdot D_1=E_1\cdot D_2=1$ and $m_1,\, m_2$ arbitrary,
\item or (up to permutation) $E_1\cdot D_1=m_1=2$ and $E_1\cdot D_2=1,\, m_2$ arbitrary.
\end{enumerate}
This last case (which is not normal crossings) is excluded by two blow-ups which make the total transform of $\Delta$ a normal crossings divisor. We indeed get the following extended dual graph (in which the $E'_i,\,i=1,2,3$ appear as white circles in the order $E'_1,\,E'_3,\,E'_2$, while the left (resp. upper) black dot is the strict transform of $D_1$ (resp. $D_2)$:
 
\begin{figure}[!h]
\setlength{\unitlength}{1.2mm}
\centering
\begin{picture}(120,18)
\put(25,1){\circle*{2}}
\multiput(35,1)(10,0){3}{\circle{2}}
\put(45,11){\circle*{2}}\put(45,2){\line(0,1){8}}
\multiput(26,1)(10,0){3}{\line(1,0){8}}
\end{picture}
\caption{A non-klt case.}
\end{figure}
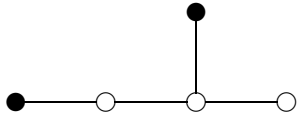

\noindent We also have: $e'_1=e_1+2,\,e'_3=1$ and $e'_2=2$. A direct resolution of the linear system $(\star)$ shows that this pair is not klt.

Assume next that our dual graph contains a fork (of white vertices). We know that $K_{X'}\equiv 0$ in this case ($e_i=2,\forall i$). Let $D_1,\,E_1,\dots,\,E_{n-2},\,E_{n-1},\,E_n$ be a shortest subgraph connecting the component $D_1$ to the fork $E_{n-2},\,E_{n-1},\,E_n$, with end points the last two vertex. The sequence is so labelled that $D_1$ meets $E_1$ only, and that $E_1,\dots,E_{n-2}$ is a chain. Consider the curve
$$B:=-(1-\frac{1}{m_1})(D_1\cdot E_1)\left(\sum_{\ell=1}^{n-2}E_{\ell}+\frac{1}{2}(E_{n-1}+E_n)\right).$$
A computation similar to the preceding one (but simpler since now $K_{X'}\equiv 0$) shows that $A\cdot E_j=B\cdot E_j,\forall j$. Thus $A=B$ and we are in the second case.

The only case left is no fork and only one component in $\Delta'$. This is the last diagram (with transversal intersection by the klt and integral conditions, as above) since the degenerate case where $n=2$ in the preceding case would lead to the existence of a fork by connectedness of $E$.
\end{demo}

\begin{rem}\label{rklt+} In the first case, the computation gives additionally: $(e_1-1)m_1=(e_n-1)m_2$ and also $e_i=2,\,\forall i\neq 1,n$. The arguments of \cite[III (5.1)]{BHPV} show that the singularity of $X$ is cyclic of type $A_{N,q}$ with $N=(n-1)(e_1-1)(e_n-1)+(e_1-1)+(e_r-1)$ and $q=(n-1)(e_1-1)+1$ with $\frac{N}{q}=(e_n-1)+\frac{e_1-1}{(n-1)(e_1-1)+1}$.
\end{rem}

\begin{lem}\label{llfg}
Let $(X,\Delta)$ be a two-dimensional germ of pair as above. This germ is klt if and only if has a finite local fundamental group. 
\end{lem}
\begin{demo}
If $\Delta=0$, it is \cite[Prop. 4.18]{KM98}. Otherwise, if the singularity is of type $A_n$, an easy computation shows (as in \cite[Ch. IV, \S 13-14]{L86}) that the local fundamental group has order $nm_1m_2$. If the singularity is of type $D_n$, we reduce it to the preceding case by a double cover.
\end{demo}

To deal with the smooth case, we need preliminary observations (\ref{loe}, \ref{eoe} and \ref{lbu} below) on orbifold \'etale covers.

\subsection{Orbifold \'etale covers.}

This first lemma ensure that orbifold \'etale covers of klt pairs will remain klt.
\begin{lem}[Prop. 5.20 \cite{KM98}]\label{loe}
Let $g:(X',\Delta')\to(X,\Delta)$ be a finite map between normal germs of varieties. Assume that $K_{X'}+\Delta'$ and $K_{X}+\Delta$ are $\QQ$-Cartier\footnote{In dimension $2$ the $\QQ$-Cartier assumption is superfluous.}, and that $K_{X'}+\Delta'=g^*(K_X+\Delta)$. Then $K_{X'}+\Delta'$ is klt if and only if $K_{X}+\Delta$ is klt.
\end{lem}

\noindent We come back in dimension $2$ with integral pairs.

\begin{defi}
We say that a finite map $g:(X',\Delta')\to(X, \Delta)$ of degree $d$ is \emph{orbifold \'etale} if
\begin{enumerate}[(i)]
\item the map $g$ ramifies only above the support of $\Delta=\sum_k (1-\frac{1}{m_k})D_k$,
\item $g$ has order of ramification $r_k$ dividing $m_k$ along $D_k$,
\item $\Delta'=g^*(\sum_k(1-\frac{r_k}{m_k})D_k)$.
\end{enumerate}
\end{defi}

\begin{rem}\label{roe}
In this preceding case the local fundamental group of $(X',\Delta')$ has index $d$ in the local fundamental group of $(X,\Delta)$. Moreover from the ramification formula we get:
$$K_X'+\Delta'=g^*(K_X+\sum_{k=1}^{r}(1-\frac{1}{r_k})D_k)+\sum_{k=1}^r(1-\frac{r_k}{m_k})g^*(\frac{D_k}{r_k})=g^*(K_X+\Delta)$$
since $(1-\frac{1}{r_k})+(1-\frac{r_k}{m_k})\frac{1}{r_k}=1-\frac{1}{m_k}$. It also follows that $(X,\Delta)$ is klt if and only if so is $(X',\Delta')$.\\
We shall apply this remark only when $r_1=m_1$, and $r_k=1$ for $r\geq 2$. In this case $\Delta'=g^*(\sum_{k=2}^{k=r}(1-\frac{1}{m_k}).D_k)$ and the component $D_1$ `disappears'.
\end{rem}

\begin{ex}\label{eoe}
We now always assume that $X$ is a smooth germ.

1. Assume that $D_1$ is a cusp of equation $x^p=y^q$ ($p,q$ coprime) and of multiplicity $m_1$ which we write symbolically as $\Delta_1=(p,q;m_1)$. Then $X'$ is the singularity of equation $z^{m_1}=y^q-x^p$. Thus if $(X,\Delta)$ is klt we have $\frac{1}{p}+\frac{1}{q}+\frac{1}{m_1}>1$. And so $(p,q,m_1)$ is up to order either $(2,2,m)$ or $(2,3, m)$, $m=3,4,5$. Moreover if the support of $\Delta$ is reducible (if $r\geq 2$) then $X'$ must be an $A_n$ singularity and the support of $\Delta'$ must have at most two components. Thus $(p,q;m_1)$ is, up to order, $(2,2; m)$ and $r\leq 3$. If $D_1$ is not {\it a priori} assumed to be a cusp, but just be given by a parametrisation $x(t)=t^q,y(t)=t^p$ with $p,q$ coprime, then since the cover of degree $m_1$ ramified exactly over $D_1$  will be a klt singularity, and so will have the usual normal form in suitable coordinates. We thus see {\it a posteriori} that $D_1$ was indeed a cusp in suitable coordinates.

2. Assume that the support of $\Delta$ has two smooth components $D_1,D_2$ tangent at order $p\geq 2$ at the origin. They have (in this order), in suitable coordinates, equations: $y=0$ and $y=x^p$. Consider the map: $g:\CC^2\to X$ given by: $g(u,v)=(x,y):=(u,v^{m_1})$. Thus using the above remark \ref{roe}, we see that making this (orbifold \'etale) cover $g$ we are lead to the case (again klt) where $$\Delta'=g^*(\sum_{k=2}^{k=r}(1-\frac{1}{m_k})D_k),$$
with $D'_2$ having equation $v^{m_1}=u^p$. In this case, letting $d$ be the $gcd$ of $(m_1,p)$, the germ $D'_2$ splits in $d$ irreducible components which are cusps of type $(p',m'_1)$ and multiplicity $m_2$ with $p':=\frac{p}{d},\, m_1':=\frac{m_1}{d}$ and are thus smooth if and only if either $p=m_1=d$ or $p=d\neq m_1$ or $m_1=d\neq p$.
\end{ex}

\subsection{The smooth case}

We consider next the remaining case in which the germ $X$ is smooth. Write $\Delta=\sum_{k=1}^{r}(1-\frac{1}{m_k})D_k$. For $k=1\dots r$, let $t_k\geq 1$ be the multiplicity of the germ $D_k$ at the origin. Thus $D_k$ is smooth if and only if $t_k=1$. We say that an irreducible germ of curve in $X\simeq \CC^2$ has a $(p,q)$-cusp at the origin if its equation in suitable coordinates is $y^q-x^p=0$ (with $p$ and $q$ coprime); in this case the multiplicity is $\inf(p,q)$.

\begin{lem}\label{lbu1}
If $(X,\Delta)$ is klt, the only possibilities for the data $r,t_k,m_k$ are the following ones.
\begin{enumerate}
\item $\forall k,\, t_k=1$. Then $r\leq 3$. If $r=3$ then $\frac{1}{m_1}+\frac{1}{m_2}+\frac{1}{m_3}>1$ (\emph{i.e.} either $(m_1,m_2,m_3)=(2,2,m_3)$ with $m_3\geq 2$ arbitrary or  $(m_1,m_2,m_3)=(2,3,m_3)$ with $2\leq m_3\leq 5$).
\item $t_k=2$ for some $k$. Then either $r=1$ and $D_k$ has a $(2,q)$-cusp at the origin with $\frac{1}{2}+\frac{1}{q}+\frac{1}{m_k}>1$; or $r=2$ (and assume $k=1$). There are two subcases:
\begin{enumerate}[(a)]
\item $D_1$ has a $(2,q)$-cusp with $m_1=2$ and $D_2$ is smooth and has intersection multiplicity $2$ with $D_1$; $m_2$ and $q$ (odd) are arbitrary.
\item $D_1$ has a $(2,3)$-cusp with $m_1=2$, $D_2$ (smooth) has intersection multiplicity $3$ with $D_1$ and $m_2=2$.
\end{enumerate}
\item $t_k=3$ for some $k$. Then $r=1$ and $D_1$ is a $(3,5)$-cusp or $(3,4)$-cusp with $m_1=2$.
\end{enumerate}
\end{lem}




\begin{demo}
Let us make a blow-up $f:X_1\to X$ at the origin with exceptional divisor $E_1$ and denote $K_1:=K_{X_1}$ and $K=K_X$. Let $\Delta_1$ be the strict transform of $\Delta$ in $X_1$. Then $K_1+\Delta_1=f^*(K+\Delta)+ cE_1$ where $c=1-\sum_k t_k(1-\frac{1}{m_k})$. The klt condition implies $c>-1$ that is $\sum_k t_k(1-\frac{1}{m_k})<2$. Since $(1-\frac{1}{m_k})\geq 1$, this implies $\sum_k t_k<4$. Thus $r\leq 3$ and we get the following list of possible values:
\begin{enumerate}[1.]
\item if $t_k\ge3$ for some $k$ then $r=1$ and $t_1=3$.
\item if $t_k=2$ for some $k$ then $r=1,2$ and $(t_1,t_2)=(2,1)$ if $r=2$.
\item if $\forall k,\,t_k=1$ then $r=1,2$ or $3$.
\end{enumerate}
We now examine the distribution of possible multiplicities.


If $t_1\geq 2$ and $r=1$ the \'etale orbifold cover (see \ref{eoe} above) of degree $m_1$ ramified along the $(p,q)$-cusp $D_1$ leads to the singularity $z^{m_1}=y^q-x^p$ with the zero orbifold divisor, which is klt if and only if the claimed inequality $\frac{1}{p}+\frac{1}{q}+\frac{1}{m_1}>1$ holds. If $t_1=2,3$, we thus get the cases described in cases 2 and 3 of lemma \ref{lbu1}. We are left with the case $t_1=2,\, t_2=1,\, r=2$. Thus $D_1$ is a $(2,q)$-cusp of multiplicity $m_1$ and $D_2$ smooth of multiplicity $m_2$ with $q$ odd. We distinguish two cases according to whether the intersection multiplicity of $D_1$ and $D_2$ is $2$ or more.

In the first case, the orbifold \'etale cover of $(X,\Delta)$ ramified to order $m_2$ along $D_2$ leads to the orbifold divisor $\Delta'$ supported on the locus of equation $x^2=y^{m_2.q}$ and with multiplicity $m_1$. From $(X',\Delta')$ still being klt we derive $\frac{1}{2}+\frac{1}{qm_2}+\frac{1}{m_1}>1$ and it yields $m_1=2$ since $qm_2\geq 6$.

If the intersection multiplicity of $D_1$ and $D_2$ is $3$ or more, the orbifold \'etale cover of $(X,\Delta)$ ramified to order $m_2$ along $D_2$ leads to the orbifold divisor $\Delta'$ supported on the locus of equation $x^{2m_2}=y^{q}$, and with multiplicity $m_1$ (still klt). If $q$ and $m_2$ are coprime we get the inequality $\frac{1}{q}+\frac{1}{2m_2}+\frac{1}{m_1}>1$ whose the only possible solution is $m_1=m_2=2$ and $q=3$.

We now show that $q$ and $m_2$ are coprime, which will complete the proof. Otherwise let $d>1$ be their $gcd$. Then $\Delta'$ consists of $d$ components of multiplicity $m_1$ and equations $x^t=\varepsilon y^s$ with $t=\frac{q}{d}$ and $s=\frac{2m_2}{d}$. These components need to be smooth since $d>1$. Thus either $t=1$ or $s=1$. Since $q$ is odd it is not divisible by $2m_2$ and so $q=d$ divides $m_2$. Because $\Delta'$ is supported by $d=q\leq 3$ smooth components, we finally get $d=q=3$ components having pairwise tangency of order $\frac{2m_2}{q}\geq 2$ at the origin. This contradicts lemma \ref{lbu} below.
\end{demo}

\begin{lem}\label{lbu}
Let $(X,\Delta)$ be a two-dimensional smooth germ, the support of $\Delta$ consisting of $s\geq 2$ smooth germs having two-by-two contact order $t\geq 1$ at the origin and each of multiplicity $m_k$ ($k=1,\dots,r$). Then $(X,\Delta)$ is klt if and only if
$$\sum_k(1-\frac{1}{m_k})<1+\frac{1}{t}.$$
The solutions $(t,r,m_1\leq m_2\leq \dots\leq m_r)$ are:
\begin{enumerate}[$\star$]
\item $t=1,\, r=2$ and $(m_1,m_2)$,
\item $t=1,\, r=3$ and $(2,2,m_3)$ or $(2,3,3\leq m_3\leq 5)$,
\item $t=2,\, r=2$ and $(2,m_2)$, $(3,m_2\leq 5)$,
\item $t=3,\, r=2$ and $(2,m_2\leq 5)$,
\item $t\geq 4,\, r=2$ and $(2,2)$.
\end{enumerate}
\end{lem}
\begin{demo}
We get a log-resolution of this pair in performing $t$ successive suitable point blow-ups with successive exceptional divisors $E_1,\dots,E_t$. If $b:X'\to X$ is the resulting composition one checks that in this process $$K_{X'}=b^*(K_X)+\sum_{h=1}^{t}hE_h$$
while
$$b^*(\Delta)=\Delta'+(\sum_k(1-\frac{1}{m_k}))(\sum_{h=1}^{t}hE_h)$$
where $\Delta'$ is the strict transform of $\Delta$. We thus see that the initial pair is klt if and only if the coefficient of $E_t$ in $K_{X'}+\Delta'-b^*(K_X+\Delta)$ is strictly greater than $-1$, \emph{i.e.} if
$$\sum_k(1-\frac{1}{m_k})<1+\frac{1}{t}.$$
The solutions are then easily determined.
\end{demo}

\begin{prop}\label{orbi lisse}
The following list exhaust the classification of integral pairs which are klt with a smooth ambient surface. We use the notation of \cite{U07}. 
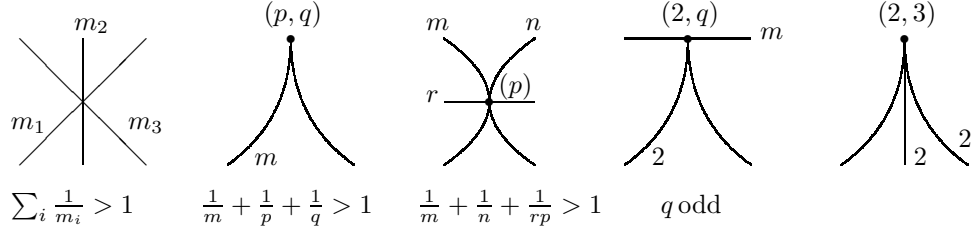
\begin{figure}[!h]
\setlength{\unitlength}{1.2mm}
\begin{picture}(120,20)
\put(1,1){\line(1,1){14}}\put(1,15){\line(1,-1){14}}\put(8,1){\line(0,1){14}}
\put(0,5){$m_1$}\put(7,16){$m_2$}\put(13,5){$m_3$}
\qbezier(24,1)(31,6)(31,15)\qbezier(31,15)(31,6)(38,1)\put(31,15){\circle*{1}}
\put(28,17){$(p,q)$}\put(27,1){$m$}
\qbezier(48,1)(58,8)(48,15)\qbezier(58,1)(48,8)(58,15)\put(53,8){\circle*{1}}\put(48,8){\line(1,0){10}}
\put(54,9){$(p)$}\put(46,8){$r$}\put(46,16){$m$}\put(57,16){$n$}
\qbezier(68,1)(75,6)(75,15)\qbezier(75,15)(75,6)(82,1)\put(75,15){\circle*{1}}\put(68,15){\line(1,0){14}}
\put(72,17){$(2,q)$}\put(83,15){$m$}\put(71,1){$2$}
\qbezier(92,1)(99,6)(99,15)\qbezier(99,15)(99,6)(106,1)\put(99,15){\circle*{1}}\put(99,1){\line(0,1){14}}
\put(96,17){$(2,3)$}\put(100,1){$2$}\put(105,3){$2$}
\end{picture}

\begin{picture}(120,5)
\put(0,1){$\sum_i \frac{1}{m_i}>1$}\put(21,1){$\frac{1}{m}+\frac{1}{p}+\frac{1}{q}>1$}
\put(45,1){$\frac{1}{m}+\frac{1}{n}+\frac{1}{rp}>1$}
\put(72,1){$q\,\mathrm{odd}$}
\end{picture}
\caption{Germs of klt pairs with smooth base (and conditions on the coefficients)}
\end{figure}
\end{prop}
\begin{demo}
To prove Proposition \ref{orbi lisse} we are thus left with the case where $X$ is smooth and $\Delta=\sum_{k=1}^{k=r}(1-\frac{1}{m_k})D_k$, the $D_k$ being smooth at the origin. We know that $r\leq 3$. When $r=1$ any $m_1$ leads to the klt situation. When $r=2$ the possible situations are described in the lemma \ref{lbu}. We thus assume that $r=3$ and that not all three components are normal crossings. 

We shall then show that (after reordering) $D_1$ and $D_2$ are tangent to order $t\geq 2$ while $D_3$ is transversal to them and that $\frac{1}{m_1}+\frac{1}{m_2}+\frac{1}{tm_3}>1$. So that one has either $m_1=m_2=2$, $m_3$ and $t$ arbitrary or $2=t=m_1=m_3$ and $m_2=3$.

Assume indeed that $D_1$ and $D_2$ have contact order $t\geq 2$. We first remark that $D_3$ is transversal to them. Assume not, then after possible reordering of the components, we may assume that the order of contact of $D_3$ with $D_1$ and $D_2$ is at least $t$. In this situation, the proof of the above lemma \ref{lbu} still applies to show that 
$\sum_{k=1}^{3}(1-\frac{1}{m_k})<1+\frac{1}{t}\leq 1+\frac{1}{2}$, which is impossible since $m_k\geq 2,\forall k$.

The conclusion now simply follows from considering the orbifold \'etale cover ramified to order $m_3$ along $D_3$, which replaces $\Delta$ by $\Delta'$, supported by two smooth components of multiplicities $m_1$ and $m_2$ and having contact order $t'=tm_3$. This implies by lemma \ref{lbu} that $\frac{1}{m_1}+\frac{1}{m_2}+\frac{1}{tm_3}>1$ since the inequality $\sum_k(1-\frac{1}{m_k})<1+\frac{1}{t'}$ is equivalent to $\sum_k\frac{1}{m_k}+\frac{1}{t'}>1$ when $r=2$. The solutions are the classical ones.
\end{demo}

\subsection{The local fundamental groups.}

\begin{prop}\label{pgfl} Let $(X,\Delta)$ be a two-dimensional germ of integral pair. If this germ is klt it has a finite local fundamental group. 
\end{prop}
\begin{demo}
If $X$ is singular this is lemma \ref{llfg}. When $X$ is smooth, this is a direct consequence of the construction of suitable orbifold \'etale covers as explained above.
\end{demo}

\bibliographystyle{amsalpha}


\end{document}